\RequirePackage{ifpdf}
\ifpdf 
\documentclass[pdftex]{sigma}
\else
\documentclass{sigma}
\fi

\begin{document}

\allowdisplaybreaks

\renewcommand{\PaperNumber}{027}

\renewcommand{\thefootnote}{$\star$}

\FirstPageHeading

\ShortArticleName{Deformation Quantization in White Noise
Analysis}

\ArticleName{Deformation Quantization in White Noise
Analysis\footnote{This paper is a contribution to the Proceedings
of the Workshop on Geometric Aspects of Integ\-rable Systems
 (July 17--19, 2006, University of Coimbra, Portugal).
The full collection is available at
\href{http://www.emis.de/journals/SIGMA/Coimbra2006.html}{http://www.emis.de/journals/SIGMA/Coimbra2006.html}}}

\Author{R\'emi L\'EANDRE}

\AuthorNameForHeading{R. L\'eandre}

\Address{Institut de Math\'ematiques, Universit\'e de Bourgogne,
21000 Dijon, France}
\Email{\href{mailto:Remi.Leandre@u-bourgogne.fr}{Remi.Leandre@u-bourgogne.fr}}

\ArticleDates{Received August 02, 2006, in f\/inal form January
30, 2007; Published online February 21, 2007}

\Abstract{We def\/ine and present an example of a deformation
quantization product on a~Hida space of test functions endowed
with a Wick product.}


\Classification{53D55; 60H40} 

\section{Introduction}
This work follows the work of Dito--L\'eandre \cite{Ditoleandre}
which was using tools of the Malliavin Calculus in order to
def\/ine a  Moyal product on a Wiener space.

Let us consider a f\/inite-dimensional symplectic manifold $M$. It
inherits from the symplectic form $\omega$ a Poisson structure
$\{\cdot,\cdot\}$ whose matrix is the inverse of the matrix of the
symplectic structure. The deformation quantization program of a
Poisson manifold was initiated by
Bayen--Flato--Fronsdal--Lichnerowicz--Sternheimer
\cite{BayenFlato1, BayenFlato2, DitoSterheimer, Maeda, Weinstein}.
The simplest case to study is the case of $\mathbb{R}^n \oplus
\mathbb{R}^{n*}$ endowed with its constant natural
 symplectic structure. This leads in particular to the notion of Moyal product.

In the case of a Hilbert space $H$, Dito \cite{Dito3} def\/ines
deformation quantization on a Hilbert space $H\oplus H^*$ endowed
with its constant canonical
 symplectic structure and he def\/ines a Moyal product on it on an appropriate algebra of functions.
The main remark is that the constant matrix of the associated
Poisson structure is still bounded.

This permits Dito--L\'eandre \cite{Ditoleandre} to def\/ine the
Moyal product on $W\oplus W^*$ where $W$ is an abstract Wiener
space. The constant symplectic form is the standard
 one on $H \oplus H^*$, the underlying Hilbert space of $W \oplus W^*$,
 such that the constant matrix of Poisson structure is still bounded.
In that work, they consider the algebra of functionals smooth in
the Malliavin sense \cite{Ikeda, Malliavin1, Malliavin2, Nualart,
Ustunel} in order to def\/ine the Moyal product on $W\oplus W^*$.

We consider in this work the case where the underlying Hilbert
space of the theory is a~Sobolev Hilbert space of maps from the
circle into $\mathbb{R}^n$ endowed
 with a constant symplectic structure. We do not consider the standard 2-form in order
 to do the quantization, but another form which is still interesting
to consider. The constant matrix of the involved Poisson structure
is unbounded so that we cannot use the construction
of~\cite{Ditoleandre}. This leads to some modif\/ications:
\begin{itemize}\itemsep=0pt
\item We replace the algebra of functionals of Malliavin type by a
Hida test functional space \cite{Berezanskii, Hida, HidaKuo,
Obata}.

\item We replace the Wiener product by the normalized Wick product
\cite{Hida}.
\end{itemize}

 The behaviour of our theory is completely dif\/ferent of the behaviour in f\/ield theory of
 the classical Garding--Wightman result \cite{Garding} for Canonical Commutation Relations
 in inf\/inite dimension. This is related to the fact that the Hida
 test algebra is  a space of continuous functionals on an abstract
 Wiener space associated to the Hilbert space of the
theory. The Hida test algebra is so small  that quantities as
white noise behaves as if we were in f\/inite dimension.
Nevertheless, there
 is a tentative to interpret some quantum
f\/ields by using tools of white noise analysis~\cite{HuangLuo,HuangRang}.

White noise analysis was created in order to understand very
singular objects as, for instance, the speed of the Brownian
motion. This explain that we get equivalences of some deformation
quantization in white noise analysis, which were inequivalent in
the theory of~\cite{Dito3}.

Our motivation comes from f\/ield theory \cite{Dito1, Dito2,
Duetsch, Witten}. We choose a simple model in this approach, but
it should be possible to consider more complicated Gaussian
models. We have chosen this simple model in order to get simple
computations. On the other hand, the free loop space is a
well-known object of conformal f\/ield theory and string theory.

\section{Deformation quantization in white noise analysis}

Let $H(S^1; \mathbb{R}^n)$ be the Hilbert space of maps $\gamma$
from the circle into $\mathbb{R}^n$ endowed with its canonical
Hilbert structure such that
\begin{equation*}
\int_0^1\vert \gamma(s)\vert^2ds + \int_0^1 \vert
\gamma'(s)\vert^2 ds = \Vert h \Vert^2 < \infty.
\end{equation*}
We get by Fourier expansion an orthonormal basis $\gamma_{k,i}$ of
this Hilbert space:
\begin{equation*}
\gamma_{k,i}(s) =\left(\sqrt{C_1k^2+1}\right)^{-1}e_i\cos(2\pi ks)
\end{equation*}
if $k \geq 0$ and if $k<0$
\begin{equation*}
\gamma_{k,i}(s) = \left(\sqrt{C_1k^2+1}\right)^{-1}e_i\sin(2\pi
ks),
\end{equation*}
where $e_i$ is the canonical basis of $\mathbb{R}^n$.

We consider a multiindex $I = ((k_1, i_1),\dots,(k_{\vert I\vert},
i_{\vert I \vert})$. We introduce the Hida weight:
\begin{equation*}
w_r(I) = \prod \left(\sqrt{C_1k_i^2+1}\right)^r
\end{equation*}
associated to this multiindex. $F^I$ denotes the normalized
symmetric tensor product of the $\gamma_{k_i,e_i}$ associated to
this multiindex.

We consider the weighted Fock space $W\cdot N_{r,C}$ of series
\begin{equation*}
\sum b_IF^I = F
\end{equation*}
such that
\begin{equation*}
\Vert F\Vert^2_{r,C} = \sum \vert b_I\vert^2w_r(I)C^{\vert I\vert}
< \infty,
 \end{equation*} where $b_I \in \mathbb{C}$.
In order to avoid some redundancies, we  order the multiindices by
lexicographic order such that after this choice $F$ is written in
a unique way.

\begin{definition}
The Hida test function $W\cdot N_{\infty-}$ space is the
intersection of all $W\cdot N_{r,C}$
 for $r>0$ and $C>0$ endowed with its natural topology.
 \end{definition}
 It is a Fr\'echet space. Since the system of norms $\Vert\cdot \Vert_{r,C}$
 increase with $r$ and $C$, it is a projective limit of Hilbert spaces.
Moreover if $r$ and $C$ are large, the imbedding of $W\cdot
N_{r,C}$ in $W\cdot N_{r',C'}$ is Hilbert--Schmidt for large $r'$
and $C'$.
 This follows from the classical fact that
\begin{equation*}
\sum w_r(I)^{-1}C^{\vert I\vert} < \infty
\end{equation*} if $r$ is large enough and $C$  small. This shows that $W\cdot N_{\infty-}$ is a nuclear space.

We def\/ine the Wick product $\colon   F^{I_1}F^{I_2} \colon$ as
the normalized symmetric tensor product of all the $\gamma_{k,i}$
in the concatenation of the multiindices $I_1$ and $I_2$.
\begin{theorem}\label{theorem1}
 $W\cdot N_{\infty-}$ is a commutative algebra for the Wick product.
 \end{theorem}
\begin{proof}
Let \begin{gather*} F_1 = \sum b_I^1F^I, \qquad
 F_2 = \sum b_I^2F^I.
 \end{gather*}
Therefore
\begin{equation*}
F_3 = \colon  F_1F_2 \colon  = \sum b_I^3F^I,
\end{equation*}
where
\begin{equation}\label{eq11}
a_I^3 = \sum_{I_1, I_2}b_{I_1}^1b_{I_2}^2,
\end{equation}
where the sum runs over all considered multiindices $I_1$ and
$I_2$ whose concatenation is $I$. If $I_1$ and $I_2$ are such
multiindices, we have clearly
\begin{equation*}
w_r(I) = w_r(I_1)w_r(I_2).
\end{equation*}
Moreover there are at most $C^{\vert I\vert} = C^{\vert
I_1\vert}C^{\vert I_2\vert}$ terms in the sum \eqref{eq11} such
that
\begin{equation*}
 \vert b_I^3\vert^2 \leq C\sum_{I_1, I_2}\big(\vert b_{I_1}^1
 \vert^2C^{\vert I_1\vert}\big)\big(\vert b_{I_2}^2\vert^2C^{\vert I_2\vert}\big).\tag*{\qed}
 \end{equation*}\renewcommand{\qed}{}
 \end{proof}

 (We refer to \cite{LeandreRogers} for an analogous statement.)

\begin{definition}A Poisson structure on $W\cdot N_{\infty-}$ $\{\cdot ,\cdot \}$
is given by a $\mathbb{C}$-bilinear map $\{\cdot,\cdot\}$ from
$WN_{\infty-}\times WN_{\infty-}$ into $W\cdot N_{\infty-}$ such
that:
\begin{enumerate}\itemsep=0pt
\item[(i)] $\{\cdot,\cdot\}$ is antisymmetric, satisf\/ies the
Jacobi derivation and is a derivation with respect of the Wick
product in each argument;

\item[(ii)] if $F_1 = 1$, $\{F_1, F_2\} = 0$; \item[(iii)] for all
$r$ and $C$, there exist $r_1$ and $C_1$ such that
\begin{equation*}
\Vert \{F_1, F_2\}\Vert_{r,C}\leq K \Vert F_1\Vert_{r_1, C_1}\Vert
F_2\Vert_{r_1, C_1}.
\end{equation*}
\end{enumerate}
\end{definition}

In this formalism we can  easily consider the  deformation
quantization of \cite{BayenFlato1, BayenFlato2, DitoSterheimer,
Maeda, Weinstein}.

We consider the set of formal series $W\cdot N_{\infty-}[[h]]$
with coef\/f\/icients in the Hida test functional space.

\begin{definition}\label{def3}
A star-product on $W\cdot N_{\infty-}[[h]]$ is a
$\mathbb{C}[[h]]$-bilinear product $*_h$ on
$WN_{\infty-}[[h]]\times W\cdot N_{\infty-}[[h]]$ with values in
$W\cdot N_{\infty-}[[h]]$ given by
\begin{equation*}
F_1*_hF_2 = \sum_{l\geq 0}h^lP_l(F_1,F_2)
\end{equation*}
for $F_1$ and $F_2$ belonging to the Hida test functional space.
The star-product is extended by $\mathbb{C}[[h]]$-bilinearity to
$W\cdot N_{\infty-}[[h]]$and satisf\/ies if $F_1$, $F_2$, $F_3$
belong to $W\cdot N_{\infty-}$:
\begin{enumerate}\itemsep=0pt
\item[(i)] $P_0(F_1,F_2) = \colon  F_1F_2 \colon$; \item[(ii)]
$P_1(F_1,F_2)-P_1(F_2,F_1) = 2\{F_1, F_2\}$; \item[(iii)] for all
$r$, $C$, $l$, there exist $r_1$, $C_1$ such that
\begin{equation*}
\Vert P_l(F_1,F_2)\Vert_{r,C}\leq K \Vert F_1\Vert_{r_1,C_1}\Vert
F_2\Vert_{r_1,C_1}
\end{equation*}
and $P_l$ vanishes on constants. \item[(iv)] $F_1*_h(F_2*_hF_3) =
(F_1*_hF_2)*_h F_3$.
\end{enumerate}
\end{definition}

\section{Example: the Hida star product}
Let $\omega = \omega_{i,j}$ be a nondegenerate antisymmetric
bilinear form on $\mathbb{R}^{n}$ ($n$ is even).
 Without loss of generality, we can write $\omega$ as:
\begin{equation}\label{eq17}
\omega_{2i,2i+1} = -\omega_{2i+1, 2i} = 1
\end{equation} and  $\omega_{i,j} = 0$ elsewhere. Namely,
 this diagonalization do not change the Hida space
 of test functionals  we have considered in this paper.
 This comes from the fact
 that a linear transformation on $\mathbb{R}^n$
 induces a linear transformation on $(\mathbb{R}^n)^{\hat{\otimes} r}$,
 the symmetric tensor product of $\mathbb{R}^n$ of length~$r$ of
 norm bounded by $C^k$.  If we perform a linear change of coordinates
 on $\mathbb{R}^n$, this induces on the system of $F_I$ $\vert I \vert = r$ a
 linear transformation.
 But it is a linear transformation  on each block realized by
 the $F^I$ $I = (k_I,i_I)$, $k_I$ being f\/ixed.
Each block is identif\/ied with $(\mathbb{R}^n)^{\hat{\otimes}
r}$, where we can look at the induced linear transformation.
 But, we could repeat the following considerations without  using this
diagonalization of $\omega$.

We introduce the antisymmetric bilinear
 form on $H(S^1;\mathbb{R}^n)$
\begin{equation*}
\Omega(\gamma^1, \gamma^2) =
\int_0^1\omega(\gamma^1(s),\gamma^2(s))ds, \qquad  \Omega =
\Omega_{(k_1,i_1), (k_2,i_2)},
\end{equation*}
where
\begin{equation*}
\Omega_{(k,2i), (k,2i+1)} =(Ck^2+1)^{-1}= -\Omega_{(k,2i+1),
(k,2i)}
\end{equation*}
and other components of $\Omega$ elsewhere vanish. This
antisymmetric bilinear form leads to a~singular Poisson structure
whose matrix form is
 $\Omega^{-1}$ which is not bounded: therefore the theory
 of~\cite{Ditoleandre} is not suitable to describe the deformation
 quantization theory related to this symplectic structure.

Let $a_{(k,i)}$ be the annihilation operator on the symmetric Fock
space associated to $\gamma_{(k,i)}$. We have that:
\begin{equation*}
a_{(k_1,i_1)}\cdots a_{(k_l,i_l)}F^I = C(I_1,
(k_1,i_1),\dots,(k_l,i_l))F^{I_1},
\end{equation*}
where in $I_1$ we have removed $(k_1,i_1),\dots,(k_l,i_l)$ if it
is possible (in the other case the expression vanishes). We have
the bound
\begin{equation}\label{eq21}
\vert C(I_1,(k_1,i_1),\dots,(k_l,i_l))\vert \leq C_l^{\vert
I\vert}.
\end{equation}

The constant $C(I_1,(k_1,i_1),\dots,(k_l,i_l))$ comes from the
fact that we are considering the Bosonic Fock space: if we
consider the same Boson at the power $n$, the associated
annihilation operator transforms it in $n$-times the Boson at the
power $n-1$.

We consider f\/inite sums $F_1 = \sum b_I^1F^I$ and $F_2 = \sum
b_I^2 F^I$. In a traditional way, we can put
\begin{equation*}
\{F_1,F_2\} = \sum\Omega^{(k_1,i_1),(k_2,i_2)}\colon
a_{(k_1,i_1)}F_1a_{(k_2,i_2)}F_2 \colon.
\end{equation*}
(We consider normalized Wick products.)
$\Omega^{(k_1,i_1),(k_2,i_2)}$ are the generic elements of the
inverse of the symplectic form $\Omega$ and are not bounded.
Therefore $\{\cdot,\cdot\}$ do not act on the Sobolev spaces of
the Malliavin Calculus  unlike the Poisson structure studied
in~\cite{Ditoleandre}, and we have to consider dif\/ferent
functional spaces if we want to extend the previous formula from
f\/inite sums  to series.   We have:

\begin{proposition}
$\{\cdot,\cdot\}$ defines a Poisson structure in the sense of
Definition~{\rm \ref{def3}} on $W\cdot N_{\infty-}$.
\end{proposition}
\begin{proof}Let us show f\/irst (iii) in Def\/inition~{\rm \ref{def3}}. We have
\begin{equation*}
\{F_1,F_2\}= \sum b_I^3F^I,
\end{equation*}
where
\begin{equation}\label{eq24}
b_I^3 = \sum
(Ck^2+1)b^1_{I_1\cup(k,2i)}b^2_{I_2\cup(k,2i+1)}C(I_1,
(k,2i))C(I_2, (k, 2i+1))+A,
\end{equation}
where $A$ is a similar term and where we sum over all $ k$, $i$,
$I_1$, $I_2$ so that the concatenation $ I_1 \cup I_2$ of $I_1$,
$I_2$ is equal to $I$.
 We apply
Cauchy--Schwartz inequality in $(k,i)$, we use the bound of
$C(I_1,(k,2i))$, $C(I_2, (k,2i+1))$ given in \eqref{eq21} to get
\begin{gather*}
\omega_r(I)\vert b_I^3\vert^2 \leq K\sum_{k,k',i,i'}
\big(\omega_{r_1}(I_1\cup(k,2i))\vert b^1_{I_1\cup(k,2i)}\vert^2C_1^{\vert I_1\vert+1}\big)\nonumber\\
\phantom{\omega_r(I)\vert b_I^3\vert^2 \leq}{}
\times\big(\omega_{r_1}(I_2\cup(k',2i'+1))\vert b^2_{I_2\cup(k',
2i'+1)}\vert^2C_1^{\vert I_2\vert+1}\big)
\end{gather*}
for  $r_1$ and  $C_1$ large enough, where we sum on the  set of
 multiindices $I_1$ and $I_2$ is such that $I_1 \cup I_2 = I$.
 We do the same for $A$ in \eqref{eq24}.
 We use for
that \begin{equation*} \sum \big(Ck^2+1\big)^{-r'} < \infty
\end{equation*} if $r'>1$.
We deduce from the previous inequality that
\begin{equation*}
\Vert  \{F_1,F_2\}\Vert_{r,C}\leq K\left(\sum \omega_{r_1}(I)\vert
I\vert C_1^{\vert I\vert}\vert b^2_I\vert^2\right) \left(\sum
\omega_{r_1}(I)\vert I\vert C_1^{\vert I \vert}\vert
b_I^2\vert^2\right).
\end{equation*}
Therefore we deduce that
\begin{equation}\label{eq28}
\Vert \{F_1,F_2\} \Vert_{r,C}\leq K \Vert F_1\Vert_{r_1, C'_1}
\Vert F_2 \Vert_{r_1,C'_1}
\end{equation}
for $r_1$, $C'_1$ large enough. This shows (iii).

The algebraic properties of the Poisson product arise from the
fact that the family of annihilation operators commute and  an
annihilation operator is a derivation for the Wick product.
\end{proof}

If $F_1$ and $F_2$ are f\/inite sums, we can def\/ine as usual by
using the Wick product:
\begin{gather*}
F_1*_hF_2 = \sum_{l\geq 0}(-h/2)^ll!^{-1}\sum\Omega^{(k_1,i_1),(k'_1,i'_1)}\cdots \Omega^{(k_l,i_l),(k'_l,i'_l)}\nonumber\\
\phantom{F_1*_hF_2 =}{}\times\colon  a_{(k_1,i_1)}\cdots
a_{(k_l,i_l)}F_1a_{(k'_1,i'_1)}\cdots a_{(k'_l,i'_l)}F_2 \colon.
\end{gather*}
The sum is in fact f\/inite since $F_1$ and $F_2$ are f\/inite
sums. It is the exponential of the Poisson Bracket.
 Let us stress the dif\/ference with~ \cite{Ditoleandre}:  in \cite{Ditoleandre}, we were considering
 the canonical symplectic form on $W\oplus W^*$ whose inverse is
 {\bf bounded}, and so the Moyal product of \cite{Ditoleandre} was acting on the {\bf big space}
 of  Malliavin test functionals. Here it is not the case.
 Let us recall the
main dif\/ference between the Malliavin test algebra and the Hida
test algebra. The Hilbert space $H(S^1;\mathbb{R}^n)$
 induces a Gaussian measure on the Banach
space $B$ of continuous functions from $S^1$ into $\mathbb{R}^n$.
The ordinary Fock space $W\cdot N_{0,1}$ coincides with the $L^2$
of this Gaussian measure. The Malliavin test algebra is
constituted of functionals almost surely def\/ined on $B$ and the
Hida test algebra is constituted of continuous functionals on $B$.
 In order to stress the dif\/ference, let us consider chaos of length $1$,
$F = \sum\limits_{\vert I \vert = 1} b_I F^I$. These chaoses of
length~1 belong to the Malliavin algebra if and only if
 $\sum \vert b_I\vert^2 < \infty$ because the $L^p$ norms
 and the $L^2$ norms are equivalent on an abstract Wiener
 space for Wiener chaos of bounded norm. It is therefore
 clear that our Poisson structure does not act on the restriction of the
Malliavin algebra constituted of chaoses of length~1.

\begin{theorem}\label{theorem2}
Formula \eqref{eq28} can be extended in a star-product in the
sense of Definition {\rm \ref{def3}}. We call it the Hida star
product
 associated to the symplectic structure given by $\Omega$.
 \end{theorem}
\begin{proof}

The algebra is the same as in the classical case \cite{Dito3}.
Only the analysis is dif\/ferent. We put
\begin{equation*}
P_l(F_1,F_2) = \sum b_I^3F^I,
\end{equation*}
where $b_I^3$ is a sum of a bounded terms of the following type:
\begin{gather*} A= \sum (Ck_1^2+1)\cdots (Ck_l^2+1) b_{I_1\cup(k_1,2i_1)\cup\cdots \cup(k_l,2i_l)}^1\nonumber\\
\phantom{A=}\times C(I_1, (k_1,2i_1),
\dots,(k_l, 2i_l))b_{I_2\cup(k'_1, 2i'_1+1)\cup\cdots\cup(k'_l, 2i'_l+1)}^2 \nonumber\\
\phantom{A=}\times C(I_2, (k'_1,2i'_1+1),\dots,(k'_l, 2i'_l+1)),
\end{gather*}
where we sum on all $k_i$, $i_l$, $i'_l$ and all multiindices
$I_1$ and $I_2$ such that their concatenation $I_1 \cup I_2$
equals $I$. By doing as before and using the estimates
\eqref{eq21}, we deduce that
\begin{gather*}
\omega_r(I)\vert A\vert^2 \leq K
\sum\big(\omega_{r_1}(I_1\cup(k_1, 2i_1)\cup\cdots \cup(k_l,
2i_l))C_1^{\vert I_1\vert+l}
\vert b^1_{I_1\cup(k_1, 2i_1)\cup\cdots\cup(k_l,2i_l)}\vert^2\big)\\
\phantom{\omega_r(I)\vert A\vert^2 \leq}{} \times
\big(\omega_{r_1}(I_2\cup(k'_1, 2i'_1+1)\cup\cdots \cup(k'_l,
2i'_l+1)) C_1^{\vert I_2\vert+l} \vert b^2_{I_2\cup(k'_1,
2i'_1)\cup\cdots \cup(k'_l,2i'_l+1)}\vert^2\big),\nonumber
\end{gather*}
where we sum on all $(k_l, 2i_l)$, all $(k_{l'}, 2i_{l'}+1)$ and
all multiindices $I_1$ and $I_2$ such that $I_1 \cup I_2 = I$.

We deduce that:
\begin{gather*}
\Vert P_l(F_1,F_2)\Vert^2_{r,C}\leq K \left(\sum
\omega_{r_1}(I)C_1^{\vert I\vert} \vert I \vert^l\ \vert
b^1_I\vert^2\right) \left(\sum \omega_{r_1}(I)C_1^{\vert I \vert}
\vert I \vert^l\vert b^2_I\vert^2\right)
\end{gather*}
for $r_1$ and $C_1$ large enough such that
\begin{equation*}
\Vert P_l(F_1,F_2)\Vert_{r,C}\leq K \Vert F_1\Vert_{r_1,C'_1}\Vert
F_2\Vert_{r_1,C'_1}
\end{equation*}
from which the result follows.
\end{proof}
\section{Equivalence of deformations}
The main dif\/ference between this work
and~\cite{Dito3,Ditoleandre} is that the space of Hida test
functionals is very small, hence
 the space of allowed deformation is very big. This implies that some
 inequivalent deformations  in the theory of~\cite{Dito3} are here equivalent.
 In order to stress the dif\/ference, we will take the model of~\cite{Ditoleandre}.

$H$ is the Hilbert space of maps from $[0,1]$ into $\mathbb{R}$
such that
\begin{equation*}
\Vert \gamma\Vert_0^2 = \int_0^1 \vert d/ds \gamma(s)\vert^2 ds <
\infty.
\end{equation*}
We consider $H \oplus H^* = H_t$ endowed with its canonical
symplectic form. We def\/ine $W_{\infty-}$ to be the space of maps
such that $\int_0^1 \vert d^r/ds^r \gamma(s)\vert^2ds < \infty$
for all $r$. It is a Fr\'echet space. We choose a convenient
Hilbert basis of $H$: if $n>0$, $\gamma_n(s) = {\sin[2\pi ns]\over
C_1n}$ and if $n<0$, $\gamma_n(s) = {\cos[2\pi ns]-1\over C_1n}$.
 Associated to this Hilbert space, and by using the convenient
 Hida weights associated to this basis, we can def\/ine
 the Hida test algebra $W\cdot N_{\infty-}$.
 We can def\/ine a Poisson structure $\{\cdot,\cdot\}$ associated to this symplectic
 form on $H_t$
 which acts continuously on $W\cdot N_{\infty-}\times W\cdot N_{\infty-}$.
 Computations are similar to the Part~III, but simpler since the matrix
 of the Poisson structure is bounded.

Let $\gamma_1 \oplus \gamma_2$ belong to $W_{\infty-}\oplus
W_{\infty-}$. We consider the Wick exponential
 $\Phi_{\gamma_1, \gamma_2}$
 \begin{equation*}
 h_1 \oplus h_2 \rightarrow \colon
 \exp[\langle h_1, \gamma_1\rangle_0 + \langle h_2,\gamma_2\rangle_0] \colon.
 \end{equation*}
 The Wick exponentials are dense in $W\cdot N_{\infty-}$.\par
 We consider an operator $A: \gamma_n \rightarrow \lambda_n \gamma_n$ with
 $\vert \lambda_n \vert \leq C |n|^\alpha$.
 According to~\cite{Dito3}, we put
 \begin{equation*}
 E_A[F_1, F_2] = \sum \colon a_i^1F_1\lambda_i a_i^2F_2 \colon
 + \sum \colon a_i^1F_2\lambda_ia_i^2F_1 \colon,
 \end{equation*}
 where $a_i^1$ are the standard annihilation operators
 in the direction of $H$ associated to $\gamma_i$ and $a_i^2$
 are standard annihilation operators in the direction of $H^* \sim H$.
 Since $|\lambda_i\vert$ are bounded by $C|i|^\alpha$
 and  we are considering the same
 Hida weight as in the f\/irst part,
 but with this new orthogonal basis,
 it follows that:
 \begin{theorem}\label{theorem3}
 $E_A$ is continuous from
 $W\cdot N_{\infty-} \times W\cdot N_{\infty-}$ into $W\cdot N_{\infty-}$.
 \end{theorem}
 We put according to~\cite{Dito3},
 \begin{equation*}
 C_1^A[F_1, F_2] = \{F_1, F_2\} + E_A[F_1,F_2]
 \end{equation*}
 and $C_r^A[F_1,F_2] = (C_1^A)^r[F_1,F_2]$ in the sense of
  bidif\/ferential operators. $C_r^A$ is still
  continuous from $W\cdot N_{\infty-}\times W\cdot N_{\infty-}$ into $W\cdot N_{\infty-}$
  (the proof is very similar to the proof of Theorem~\ref{theorem1}).
 \begin{definition}

 We put
 \begin{equation*}
 F_1*_h^AF_2 = \colon FG \colon
  + \sum_{r\geq 1} {h^r\over r!} C_r^A(F_1,F_2).
  \end{equation*}
  \end{definition}
 As in Theorem~\ref{theorem2}, due to the polynomial growth of the $\lambda_i$, $*_h^A$
 def\/ines a quantization by deformation in Hida sense of $\{\cdot,\cdot\}$.
 But unlike in~\cite{Dito3}, we have:
 \begin{proposition} $*_h^A$ and $*_h$ are equivalent
 on the Hida test functional space.
 \end{proposition}
 \begin{proof} We put as in~\cite{Dito3}
 \begin{equation*}
 T_1F = - \sum \lambda_ia_i^1a_i^2F.
 \end{equation*}
 Due to the polynomial growth of the $\lambda_i$,
 $T' = \exp[hT_1]$ is continuous on $W\cdot N_{\infty-}[[h]]$
 Moreover, let us recall that by~\cite{Dito3} formula \eqref{eq17}
 \begin{gather*} \Phi_{\gamma_1,\gamma_2}*_h^A\Phi_{\gamma'_1, \gamma'_2} =
  \exp[h(\langle \gamma'_2, (A+\mathbb{I})\gamma_1\rangle_0
  +\rangle\gamma_2, (A+\mathbb{I})\gamma'_1\rangle_0)]\Phi_{\gamma_1+\gamma'_1, \gamma_2+\gamma'_2}.
  \end{gather*}
 We conclude as in~\cite{Dito3}, by remarking that
 \begin{equation*} T'(\Phi_{\gamma_1,\gamma_2}*_h^A
 \Phi_{\gamma'_1, \gamma'_2})
 = T'\Phi_{\gamma_1, \gamma_2}*_hT'\Phi_{\gamma'_1,\gamma'_2}.
 \end{equation*}
 This proves the theorem since the Wick exponentials are dense in the Hida space.
 \end{proof}

\begin{remark}
If $A = \mathbb{I}$, we get the normal product. Let us stress the
dif\/ference with the theory of~\cite{Dito3}. In~\cite{Dito3},
$*_h^A$ and $*_h$
 were equivalent if and only if $A$ is a Hilbert--Schmidt
operator, then the Moyal product and the normal product ($A=
\mathbb{I}$) were inequivalent. In the case of $\mathbb{R}^n$, the
Moyal product and the normal product are equivalent. For the Hida
Calculus, deformation theory behaves more or less as in f\/inite
dimension. This comes from
 that the Hida test functional space is so small that all
 algebraic considerations in f\/inite dimension, where we were considering f\/inite sums,
 remain true in this context.
\end{remark}

\pdfbookmark[1]{References}{ref}
\LastPageEnding

\end{document}